
\documentstyle{amsppt}
\topmatter
\author Piotr Kosi\' nski \footnote[*]{supported by KBN grant
no 2P30221706p02\hfill{ }},
Micha\l \ Majewski$^*$,
Pawe\l \
Ma\' slanka \footnote[**]{supported by University
of \L \' od\' z grant no 505/444\hfill{ }}
\endauthor
\title Representations of generalized oscillator algebra \endtitle
\address
Piotr Kosi\' nski, Micha\l \ Majewski \newline
Department of Theoretical Physics \newline
University of \L\' od\' z \newline
Pomorska 149/153 \newline 90-236 \L\' od\' z, Poland \newline
\newline
Pawe\l \ Ma\' slanka \newline
Institute of Mathematics \newline
University of \L \' od\' z \newline
Banacha 22 \newline
90-238 \L \' od\' z, Poland
\endaddress
\abstract
The representations of the oscillator algebra introduced by
Brzezi\' nski et al. (Phys.~Lett.~{\bf B~311} (1993), 202) are classified.
\endabstract
\endtopmatter
\TagsOnRight
\NoRunningHeads
\NoBlackBoxes
\document

\head I.~Introduction \endhead

In this paper we present an exhaustive discussion of representations
of the oscillator algebra introduced in Ref.~[1]. This algebra was
obtained as a $q$-deformation of the modified oscillator algebra
underlying the structure of Calogero model [2]. In the paper [1]
the Fock representations of the algebra under consideration
were found. However, it is known [3] that $q$-deformed oscillator
algebra posses some exotic representations which disappear in the limit
$q\to 1$. Therefore, it is expected that some phenomenon occurs for
the algebra considered by Brzezi\' nski et al. Following the methods
used in Ref.~[3] we show that this is indeed the case and provide the
complete classification of its representations.

\head II.~ Construction of representations \endhead

The algebra under consideration reads~:
$$\aligned
& aa^+-qa^+a=q^{-N}(1+2\alpha K), \\
& [N,a]=-a, \quad [N,a^+]=a^+,  \\
& \{ K,a\}=0, \quad \{ K,a^+ \} =0, \\
& [N,K]=0, \\
& N^+=N, \quad K^+=K,
\endaligned \tag1$$
where $q\in \Bbb R_+$, $\alpha \in \Bbb R -\{ 0\}$. This algebra
posses the following Casimir operators
$$C_1=K^2,\qquad C_2=Ke^{i\pi N}, \qquad C_3=e^{2i\pi N}.\tag2$$
Obviously, they are not independent
$$C_1C_3=C_2^2.\tag3$$

We will be looking for irreducible representations. Let $\gamma$
be an eigenvalue of $C_2$ corresponding to a given representation;
then
$$K=\gamma e^{-i\pi N}.\tag4$$
Let $\Psi_0$ be a common eigenvector of $N$ and $K$~:
$$\aligned
& N\Psi_0=\nu_0 \Psi_0,\\
& K\Psi_0=\gamma e^{-i\pi\nu_0}\Psi_0.
\endaligned \tag5$$
Due to the commutativity of $a^+a$ and $aa^+$ with $N$ and $K$
we may assume that
$$\aligned
& a^+a\Psi_0=\lambda_0 \Psi_0, \\
& aa^+ \Psi_0=\mu_0 \Psi_0
\endaligned \tag6$$
and $(\Psi_0,\Psi_0)=1$. It is easy to see that the vectors
$\Phi_n$ defined by
$$\Phi_n=\cases
\left( a^+ \right)^n \Psi_0 & \quad \text{for} \quad  n\geq 0, \\
a^{-n}\Psi_0 & \quad \text{for} \quad n<0
\endcases$$
are eigenvectors of $a^+a$ and $aa^+$~:
$$\aligned
& a^+a \Phi_n=\lambda_n \Phi_n, \\
& aa^+ \Phi_n=\mu_n \Phi_n.
\endaligned \tag7$$
Now, let us define the following vectors~:
$$\Psi_n=
\cases
\frac{\dsize 1}{\sqrt{\dsize{\prod_{k=1}^n \lambda_k}}}\left( a^+ \right)^n
\Psi_0\quad \text{for} \quad n\ge 0, \\
\frac{\dsize 1}{\sqrt{\dsize{\prod_{k=1}^{-n} \lambda_{n+k}}}} a^{-n} \Psi_0
\quad \text{for} \quad n<0.
\endcases \tag8$$
They are orthogonal (as eigenvectors of $N$ corresponding to
different eigenvalues) and normalized. The action of basic
operators is given by
$$\aligned
& a^+\Psi_n =\sqrt{\lambda_{n+1}}\Psi_{n+1}, \\
& a\Psi_n =\sqrt{\lambda_n}\Psi_{n-1}, \\
& N\Psi_n =(\nu_0+n) \Psi_n, \\
& K\Psi=\frac{(-1)^n}{2\alpha }B\Psi_n,
\endaligned \tag9$$
where, for later convenience, we have defined
$$B=2\alpha \gamma e^{-i\pi \nu_0} \in \Bbb R .\tag10$$

The only additional condition we have to take into account
is that $\lambda_n$ and $\mu_n$, being eigenvalues
of nonnegative operators, should be nonnegative.
Using the basic commutation rules applied to $\Psi_n$
we obtain
$$\mu_n-q\lambda_n =-q^{-(n+\nu_0)}
\left( {1+2\alpha \gamma e^{-i\pi (n+\nu_0)}} \right)  .\tag11$$
But $a\left( a^+a \right) \Psi_n =\left( aa^+ \right) a\Psi_n$
which gives
$$\lambda_n =\mu_{n-1}. \tag12$$
Equations (11) and (12) imply the following recurrence relation
$$\lambda_{n+1}=q\lambda_n +q^{-\nu_0 -n}
\left( 1+(-1)^nB \right) ,\tag14$$
which can be explicitly solved to yield~:
$$\lambda_n=\lambda_0 q^n +q^{-\nu_0}
\left( \frac{q^n-q^{-n}}{q-q^{-1}}+
B \frac{q^n-(-1)^{n}q^{-n}}{q+q^{-1}} \right) .\tag15$$
Nonnegativity of $\lambda_n$ implies
$$\lambda_0 q^{\nu_0} + \frac{1}{q-q^{-1}} +\frac{B}{q+q^{-1}}
\geq
q^{-4k}\left( \frac{1}{q-q^{-1}}+\frac{B}{q+q^{-1}} \right) ,
\tag16a$$
$$\lambda_0 q^{\nu_0} + \frac{1}{q-q^{-1}} +\frac{B}{q+q^{-1}}
\geq
q^{-(4k+2)}\left( \frac{1}{q-q^{-1}}-\frac{B}{q+q^{-1}} \right) .
\tag16b$$
We have now to distinguish several cases.
\roster
\item"{(i)}" Assume $q>1$. Then at least one of the numbers
$\dsize \frac{1}{q-q^{-1}}\pm \frac{B}{q+q^{-1}}$ is positive.
Therefore, there exists $n_0$ such that for even and/or odd
\linebreak
$n<n_0$, $\lambda_n <0$, which implies $a\Psi_n=0$ for some
$n\leq n_0$. After possible renumbering we may assume
$$a\Psi_0=0,\qquad \lambda_0=0 \tag17$$
Therefore, the representation is spanned by the vectors
$\Psi_n$, $n\ge 0$, and $\lambda_n$ are given by
$$\lambda_n =q^{-\nu_0 +n}
\left( \frac{1-q^{-2n}}{q-q^{-1}}
+B\frac{1-(-1)^n q^{-2n}}{q+q^{-1}} \right) \tag18$$
The condition $\lambda\ge 0$ gives the following restriction
on the possible values of $B$~:
$$B\ge -1 \tag19$$
However, $B=-1$ must be considered separately.
In this case $\lambda_1=0$, i.e. $\mu_0=0$ which,
together with $\lambda_0=0$ and irreducibility implies
$$a=a^+=0, \qquad N=\nu_0, \qquad K=-\frac{1}{2\alpha}.
 \tag20$$
This representation is onedimensional. For $B>-1$
the representation is spanned by the vectors
$\{ \Psi_n \}_{n=0}^{\infty}$. We shall call it the Fock
representation for obvious reasons. It is given by
equations (9) and (18) with $n\ge 0$.
\item"{(ii)}" $q<1$ and one (and only one) of the values
$\dsize \frac{1}{q-q^{-1}}\pm \frac{B}{q+q^{-1}}$ is positive.
In this case there exists $n_0$ such that for $n> n_0$ $\lambda_n$
is negative for even or odd $n$. This implies $a^+\Psi_n=0$ for
some $n\geq n_0$. After possible renumbering we get
$$a^+\Psi_0 =0 .\tag21$$
In order to find the restrictions on possible values of $B$
we note that the condition (21) implies $\mu_0$, i.e.
$\lambda_1=0$ or
$$\lambda_0=-q^{-\nu_0 -1}(1+B) ,\tag22$$
which gives $B\le -1$. For $B=-1$ we obtain
onedimensional representation (20). If $B<-1$ we have to
consider the restrictions on $B$ following from the formula
$$\lambda_n=q^{n-\nu_0} \left( -q^{-1}(1+B)
+\frac{1-q^{-2n}}{q-q^{-1}}
+B\frac{1-(-1)^n q^{-2n}}{q+q^{-1}} \right) .\tag23$$
The condition $\lambda_n \ge 0$ implies
$\dsize B\le \frac{q+q^{-1}}{q-q^{-1}}$. For
$$B< \frac{q+q^{-1}}{q-q^{-1}}\tag24$$
we have $\lambda_n>0$ and the representation is given
by equations (9) and (23) with $n\le 0$. We call this
representation anti-Fock one. For
$$B= \frac{q+q^{-1}}{q-q^{-1}}\tag25$$
all $\lambda_n$ with odd $n$ are zero. Therefore the representation
is twodimensional and given by
$$\aligned
a\Psi_0=\sqrt{\frac{2q^{-\nu_0}}{q^{-1}-q}}\Psi_{-1},
  \qquad & a^+\Psi_0=0, \\
a^+\Psi_{-1}=\sqrt{\frac{2q^{-\nu_0}}{q^{-1}-q}}\Psi_0,
  \qquad & a\Psi_{-1}=0, \\
N\Psi_0=\nu_0 \Psi_0, \qquad & N\Psi_{-1}=(\nu_0-1)\Psi_{-1}, \\
K\Psi_0=\frac{q+q^{-1}}{2\alpha \left( q-q^{-1} \right)}, \qquad
&  K\Psi_{-1}=\frac{q+q^{-1}}{2\alpha \left( q-q^{-1} \right)}\Psi_{-1}.
\endaligned \tag26 $$
\item"{(iii)}" $q<1$ and both values
$\dsize \frac{1}{q-q^{-1}}\pm \frac{B}{q+q^{-1}}$ are nonpositive
(at least one must be strictly negative). There are now the following
possibilities~:
\item"{(a)}" $$\dsize \lambda_0 q^{\nu_0}
+\frac{1}{q-q^{-1}}+ \frac{B}{q+q^{-1}}<0\tag27$$
There exists then $n_0$ such that $\lambda_n<0$ for
$n<n_0$, $n$ even or odd. Therefore the representation is given
by formulae (9), (18); it is a Fock one. To provide $\lambda_n\ge 0$
for $n\ge 0$ we have to restrict $B$ to lie in the interval
$$-1\le B<-\frac{q+q^{-1}}{q-q^{-1}} .\tag28$$
For $B=-1$ we get again onedimensional representation (20).
\item"{(b)}" $$\dsize \lambda_0 q^{\nu_0}
+\frac{1}{q-q^{-1}}+ \frac{B}{q+q^{-1}}>0\tag29$$
Equation (29) implies $\lambda_n>0$ for all $n\in \Bbb Z$.
The representation is given by equations (9), (15) with
$n\in \Bbb Z$.
\item"{(c)}" $$\dsize \lambda_0 q^{\nu_0}
+\frac{1}{q-q^{-1}}+ \frac{B}{q+q^{-1}}=0\tag30$$
If $\dsize | B|<-\frac{q+q^{-1}}{q-q^{-1}}$ all $\lambda_n>0$
and representation has the same form as in point (b).
For $\dsize B=-\frac{q+q^{-1}}{q-q^{-1}}$ all $\lambda_n$ with
$n$ even are vanishing; therefore, the representation is
twodimensional and given by formulae~:
$$\aligned
a^+\Psi_0=\sqrt{\frac{2q^{-\nu_0-1}}{q^{-1}-q}}\Psi_1,
  \qquad & a\Psi_0=0, \\
a^+\Psi_1=0, \qquad &
a\Psi_1=\sqrt{\frac{2q^{-\nu_0-1}}{q^{-1}-q}}\Psi_0, \\
N\Psi_0=\nu_0 \Psi_0, \qquad & N\Psi_1=(\nu_0+1)\Psi_1, \\
K\Psi_0=-\frac{q+q^{-1}}{2\alpha \left( q-q^{-1} \right)}\Psi_0, \qquad
&  K\Psi_1=\frac{q+q^{-1}}{2\alpha \left( q-q^{-1} \right)}\Psi_1.
\endaligned \tag31 $$
For $\dsize B=\frac{q+q^{-1}}{q-q^{-1}}$ all $\lambda_n$ with $n$
odd vanish. The representation is twodimensional and is given
by formulae (26).
\endroster

\head III.~Discussion \endhead

Let us summarize the results obtained in section II.
For $q>1$ there are two possibilities~:
\roster
\item"{(a)}" if $B>-1$ the spectrum of $N$ is bounded from below
-- we have Fock representation which is irreducible and determined by
the choice of $\nu_0$ and $B$; different choices correspond to
inequivalent representations,
\item"{(b)}" if $B=-1$ we get onedimensional irreducible
representation labeled by the values of $\nu_0$; again different
values of $\nu_0$ correspond to inequivalent representations.
\endroster
The case $q<1$ is more involved. The following possibilities have to be
distinguished~:
\roster
\item"{(a)}" for
$\dsize B<\frac{q+q^{-1}}{q-q^{-1}}$
we get irreducible anti-Fock representation; the representations are
labelled by pairs $(\nu_0,B)$ and different choices correspond
to inequivalent representations,
\item"{(b)}" for
$\dsize B=\frac{q+q^{-1}}{q-q^{-1}}$
one obtains twodimensional representations parametrized by $\nu_0$;
different values of $\nu_0$ correspond to inequivalent representations,
\item"{(c)}" for $B=-1$ we obtain again onedimensional representation
parametrized by $\nu_0$; it has the same form as for $q>1$,
\item"{(d)}" for
$\dsize -1<B<-\frac{q+q^{-1}}{q-q^{-1}}$
the representations are the Fock ones parametrized by $\nu_0$ and $B$;
for different values of these parameters we obtain inequivalent
representations,
\item"{(e)}" for
$\dsize B=-\frac{q+q^{-1}}{q-q^{-1}}$
the representations, parametrized by $\nu_0$, are twodimensional
and mutually inequivalent,
\item"{(f)}" finally, there exists a set of infinitedimensional
representations for which the spectrum of $N$ extends infinitely
in both directions. They correspond to
$\dsize |B|<-\frac{q+q^{-1}}{q-q^{-1}}$
and
$\dsize \lambda_0 q^{\nu_0}+\frac{1}{q-q^{-1}}
+\frac{B}{q+q^{-1}}\ge 0$
or
$\dsize |B|=-\frac{q+q^{-1}}{q-q^{-1}}$
and
$\dsize \lambda_0 q^{\nu_0}+\frac{1}{q-q^{-1}}
+\frac{B}{q+q^{-1}}> 0$.
It is easy to see that two such representations, labeled by
$(\nu_0,B,\lambda_0)$ and $(\nu_0',B',\lambda_0')$,
are equivalent iff
$\nu_0'=\nu_0+n$, $B'=(-1)^nB$,
$\dsize \lambda_0'=
\lambda_0 q^n+q^{-\nu_0}
\left( \frac{q^n-q^{-n}}{q-q^{-1}}
+B\frac{q^n-(-1)^n q^{-n}}{q+q^{-1}} \right) $
for some integer $n$.
\endroster

Let us now consider the limits $q\to 1$ or $B\to 0$.
It is easy to see that only the onedimensional and Fock representations
survive the limit $q\to 1$. On the other hand the $B\to 0$ limit
coincides with the results obtained in Ref.~[3]. Finally, the limit
$q\to 1$, $B\to 0$ leaves only Fock representation as it should be.
The results obtained are summarized in table 1.

\newpage
\head Table 1. \endhead

$$\matrix
\matrix \bold{Type \quad of} \\ \bold{representation} \endmatrix
& q
& \matrix \bold{Restrictions} \\ \bold{on }\quad B \endmatrix
& \matrix \bold{Restrictions} \\ \bold{on }\quad \lambda_0
                         \quad \bold{ and }\quad \nu_0 \endmatrix
& \matrix q\to 1 \\ \bold{limit} \endmatrix
& \matrix \alpha \to 0 \\ \bold{limit} \endmatrix
\\
\\
\\
\\
\\
\text{onedimensional}
& \text{arbitrary}
& B=-1
& \matrix \lambda_0=0\\ \nu_0\text{ --arbitrary}\endmatrix
& \text{exists}
& \matrix \text{does not} \\ \text{exist} \endmatrix
\\
\\
\\
\\
\\
\text{twodimensional}
& q<1
& \matrix B=\frac{q+q^{-1}}{q-q^{-1}} \\ \\ \\
         B=-\frac{q+q^{-1}}{q-q^{-1}} \endmatrix
& \matrix
     \matrix \lambda_0=\frac{2q^{-\nu_0}}{q^{-1}-q} \\
             \nu_0\text{ --arbitrary}
     \endmatrix
     \\ \\ \\
     \matrix \lambda_0=0 \\
             \nu_0\text{ --arbitrary}
     \endmatrix
  \endmatrix
& \matrix
     \matrix \text{does not} \\
             \text{exist}
     \endmatrix
     \\ \\ \\
     \matrix \text{does not} \\
             \text{exist}
     \endmatrix
  \endmatrix
& \matrix
     \matrix \text{does not} \\
             \text{exist}
     \endmatrix
     \\ \\ \\
     \matrix \text{does not} \\
             \text{exist}
     \endmatrix
  \endmatrix
\\
\\
\\
\\
\\
\text{Fock}
& \matrix q>1 \\ \\ \\ \\  q<1 \endmatrix
& \matrix B>-1 \\ \\ \\ \\
          -\frac{q+q^{-1}}{q-q^{-1}}>B>-1 \endmatrix
& \matrix
     \matrix \lambda_0=0 \\
             \nu_0\text{ --arbitrary}
     \endmatrix
     \\ \\ \\
     \matrix \lambda_0=0 \\
             \nu_0\text{ --arbitrary}
     \endmatrix
   \endmatrix
& \matrix \text{exists} \\ \\ \\ \\ \text{exists} \endmatrix
& \matrix \text{exists} \\ \\ \\ \\ \text{exists} \endmatrix
\\
\\
\\
\\
\\
\text{anti-Fock}
& q<1
& B<\frac{q+q^{-1}}{q-q^{-1}}
& \matrix \lambda_0=-q^{\nu_0 -1}(1+B) \\
          \nu_0 \text{ --arbitrary}
  \endmatrix
& \matrix \text{does not} \\ \text{exist} \endmatrix
& \matrix \text{does not} \\ \text{exist} \endmatrix
\\
\\
\\
\\
\\
\matrix \text{unbounded} \\ \text{in both}
        \ \text{directions} \endmatrix
& q<1
& \matrix |B|<-\frac{q+q^{-1}}{q-q^{-1}} \\ \\ \\
          |B|=-\frac{q+q^{-1}}{q-q^{-1}}
  \endmatrix
& \matrix
\lambda_0q^{\nu_0} +\frac{1}{q-q^{-1}}+\frac{B}{q+q^{-1}} \ge 0 \\ \\ \\
\lambda_0q^{\nu_0} +\frac{1}{q-q^{-1}}+\frac{B}{q+q^{-1}} >0
  \endmatrix
& \matrix
     \matrix \text{does not} \\
             \text{exist}
     \endmatrix
     \\ \\ \\
     \matrix \text{does not} \\
             \text{exist}
     \endmatrix
  \endmatrix
& \matrix
     \text{exists}
          \\ \\ \\
     \matrix \text{does not} \\
             \text{exist}
     \endmatrix
  \endmatrix
\endmatrix$$

\newpage
{\tensmc Remark} to the case of representations unbounded
in both directions~: they are equivalent
$$(\nu_0',\lambda_0',B')\sim (\nu_0,\lambda_0,B)$$
\centerline{iff}
$$\aligned
\dsize \nu_0'=\nu_0 & +n,\quad B'=(-1)^nB, \\
\lambda_0'  =\lambda_0q^n
+ q^{-\nu_0} & \left( \frac{q^n-q^{-n}}{q-q^{-1}}
+ B\frac{q^n-(-1)^nq^{-n}}{q+q^{-1}} \right) , \\
 & n\in \Bbb Z .
\endaligned$$


\Refs
\ref \key 1 \by T.~Brzezi\' nski, J.~L.~Egusquiza, A.~J.~Macfarlane
\jour Phys. Lett. \vol B 311 \yr 1993 \pages 202
\endref
\ref \key 2 \by L.~Brink, T.~H.~Hansson, M.~A.~Vasiliev
\jour Phys. Lett. \vol B 286 \yr 1992 \pages 109
\endref
\ref \key 3 \by M.~A.~Vasiliev \jour Int. J. Mod. Phys.
\vol A 6 \yr 1991 \pages 1115
\endref
\ref \key 4 \by G.~Rideau \jour Lett. Math. Phys.
\vol 24 \yr 147 \pages 147
\endref
\endRefs
\enddocument